\documentclass{procmult}
\usepackage{a4}
\usepackage{epsfig}
\usepackage{amssymb}
\usepackage{amsbsy}
\usepackage{amsmath}
\usepackage{amsfonts}
\usepackage{graphicx}

\title{An Approach Towards the Proof of the Strong Goldbach's Conjecture for Sufficiently Large Even Integers
}

\author{Ahmad Sabihi}

\institute{\normalsize\itshape Teaching professor and researcher at some universities of Iran}

\begin{document}

\maketitle

\begin{abstract}
We approach a new proof of the strong Goldbach's
conjecture for sufficiently large even integers by applying the
Dirichlet's series. Using the Perron formula and the Residue Theorem in
complex variable integration, one could show that any large even
integer is demonstrated as a sum of two primes. In this paper,the
Riemann Hypothesis is assumed to be true in throughout the paper.
A novel function is defined on the natural numbers set.This
function is a typical sieve function.Then based on this
function,several new functions are represented and using the Prime
Number Theorem,Sabihi's theorem, and the Sabihi's second
conjecture,the strong Goldbach's conjecture is proved for
sufficiently large even integers.
\end{abstract}
\section{INTRODUCTION}
The strong Goldbach's conjecture refers to Christian Goldbach's
theory which was represented by him with a letter to  the great
Swiss mathematician Leonard Euler in 1742. He asked him for
solving this problem. He spent much time trying to prove it, but
never succeeded. In this theory, he claims that "every even
integer greater than two could be a sum of two prime numbers". The
word "strong" is opposition to word "weak" where refers to "every
odd number greater than five could be a sum of three prime
numbers". Since two hundreds years ago, several relevant proofs on
this subject have been conducted. Using the sieve method, a lot of
mathematicians as Brun [1] in 1920 have verified various results
who expressed � every even integer is a sum of product at most 9
primes and product at most another 9 primes". Wang [2] expressed
in 1962, "every even integer is a sum of 1 prime and product at
most 4 primes". Pan [3] obtained in 1962 that every even integer
is a sum of 1 prime and product at most 5 primes. Richert [4]
obtained same result in 1969, but product at most 3 primes instead
at most 5 primes. Chen [5]in 1973 proved same for 1 prime and
product at most 2 primes. R'enyi [6] in 1947 represented that
every even integer is a sum of 1 prime and product at most �c�
primes. Montgomery and Vaughan found an exceptional set in
Goldbach's problem by the circle method [7]. We presented [8-9]
by the two papers in 2009 and 2010 stating 5 new Functions,11
new Lemmas, a novel explicit formula for computing $\pi(x)$,
several methods, a new conjecture and some new experimental
computations to prove Goldbach�s conjecture experimentaly. C.D.
Pan and C.B. Pan in [10] made the proof of some Theorems in
analytic number Theory so that we make use of some of them here. Most of these
Theorems are about the proof of the Prime Number Theorem.
\setcounter{equation}{0}
\section{ELEMENTARY DEFINITIONS}
The following new function is defined as $n$ be a natural number
and $N$ an even positive integer,then
\begin{equation}
\eta(n)=\left\{
\begin{array}{ll}
0~~~ & \mbox{if~~~~$N\equiv n~mod(p)$~~$and~~~p\leq
\sqrt{N}$}\\\\1~~~ & \mbox{if~~~~$otherwise$}\end{array} \right.\
\end{equation}\\
where $p$ denotes a prime number and $gcd(p,N)=1$.From the
function $\eta(n)$ is recognized properties of sieve function.
Just, based on the above function and known functions as
$\zeta(s)$,$\theta(x)$,$\pi(x)$ and $\psi(x)$ one is able to
define the several new following functions (2-2) to (6-2).
\begin{equation}
\zeta(s,\eta)=\sum_{n=1}^{\infty}\frac{\eta(n)}{n^{s}}\
\end{equation}
\begin{equation}
\theta(x,\eta)=\sum_{p\leq x}\eta(p)log(p)\
\end{equation}
\begin{equation}
\pi(x,\eta)=\sum_{p\leq x}\eta(p)\
\end{equation}
\begin{equation}
\psi(x,\eta)=\sum_{n\leq x}\eta(n)\Lambda(n)\
\end{equation}
\begin{equation}
\Psi(s,\eta)=\sum_{n=1}^{\infty}\frac{\eta(n)\Lambda(n)}{n^{s}}\
\end{equation}
where
\begin{equation}
\Lambda(n)=\left\{
\begin{array}{ll}
log(p)~~~ & \mbox{if~~~~$log(n)=mlog(p)$}\\\\0 &
\mbox{if~~~$otherwise$}\end{array} \right.\
\end{equation}
\begin{equation}
\theta(x)=\sum_{p\leq x}log(p)\
\end{equation}
\begin{equation}
\pi(x)=\sum_{p\leq x}1\
\end{equation}
\begin{equation}
\psi(x)=\sum_{n\leq x}\Lambda(n)\
\end{equation}
The main scope of this paper is to show that $\pi(N,\eta)\geq 1$ where
$N$ denotes a sufficient large even integer.This proves
the strong Goldbach's conjecture for sufficiently even
integers.\\\\
\textbf{Lemma 2.1}\\\\ \textit{If $\pi(N,\eta)\geq 1$, then the strong Goldbach's conjecture holds}\\\\
\textit{\textbf{Proof:}}\\\\ let
\begin{equation} \pi(N,\eta)=\sum_{p\leq N}\eta(q)=\sum_{p\leq N}\left\{
\begin{array}{ll}
0 & \mbox{if~~~~$N\equiv q~mod(p)$~~$and~~~p\leq \sqrt{N}$}\\\\1 &
\mbox{if~~~~$otherwise$}\end{array} \right.\
\end{equation}\\
Consider a contradiction as:$\pi(N,\eta)=0$\\
~~~~If $\pi(N,\eta)=0$ then consider a typical prime as $2<q<N$,
therefore based on the above relation, one should have for any q:
$N-q<\sqrt{N}$ or $N-q>\sqrt{N}$.If $N-q<\sqrt{N}$ then $N-q=kp$
and based on the equation (1-2), one concludes that $\eta(q)=0$.If
$N-q>\sqrt{N}$, one has two cases:the first $N-q$ is prime and
based upon the equation (1-2), $\eta(q)=1$. In the second
case,$N-q$ is not a prime but at least one of its factors will be
under $\sqrt{N}$(every integer as $m$ has at least a prime factor
less than $sqrt(m)$, because otherwise product of its prime
factors will be greater than itself i.e.$m$). Hence, based on the
equation (1-2), $\eta(q)=0$.From this argument, one could conclude
if $\pi(N,\eta)\geq 1$ then the first case will be happened. This
means that $N-q$ is a prime number and the strong Goldbach's
conjecture holds. Using the Prime Number Theorem [10] could be
written:
\begin{equation}
\theta(x,\eta)=\pi(x,\eta)log(x)-\int_{2}^{x}\frac{\pi(u,\eta)}{u}du\
\end{equation}
\begin{equation}
\pi(x,\eta)=\frac{\theta(x,\eta)}{log(x)}+\int_{2}^{x}\frac{\theta(u,\eta)}{u(log(u))^{2}}du\
\end{equation}
\begin{equation}
\psi(x,\eta)=\theta(x,\eta)+O(x^{1/2})\
\end{equation} \setcounter{equation}{0}\
\section{LEMMAS AND THEOREMS TO APPROACH TO THE PROOF OF
$\pi(N,\eta)\geq 1$} In this section, we prove four Lemmas
using both Perron formula of Dirichlet�s series and the Residue
Theorem. These Lemmas are the base of a basic Theorem which will
be given in the next Section. Here, we give a theorem and
some lemmas along with their proofs.\\\\
\textbf{Lemma 3.1}\\\\ \textit{Let
$A(s)=\sum_{n=1}^{\infty}a(n)n^{-s}$,$\sigma_{a}<+\infty$ and 
$H(u)$ and $B(u)$ be increasing functions so that $|a(n)|\leq
H(n)$,$n=1,2,...$ and $\sum_{n=1}^{\infty}a(n)n^{-\sigma}\leq
B(\sigma)$,$\sigma>\sigma_{a}$ for any $s_{0}=\sigma_{0}+it_{0}$
and $b_{0}\geq \sigma_{0}+b>\sigma_{a}$,$T\geq 1$. Also,$x\geq 1$,
then we have Perron formula as follows [10]:}\
\begin{eqnarray}
\sum_{n\leq
x}a(n)n^{-s_{0}}+\frac{1}{2}a(x)x^{-s_{0}}=\frac{1}{2\pi
i}\int_{b-iT}^{b+iT}A(s_{0}+s)\frac{x^{s}}{s}ds+O(\frac{x^{b}B(b+\sigma_{0})}{T})+\nonumber\\O(x^{1-\sigma_{0}}H(2x)min(1,\frac{log(x)}{T}))~~~~~~~~~~~~~~~~~~~~~~~~~~~~~\
\end{eqnarray}
\textit{where "O" denotes a constant depending on $\sigma_{a}$ and
$b_{0}$. $\sigma_{a}$ denotes absolute convergence abscissa of
Dirichlet's series.} Refer to [10] in order to see the proof.\\\\
If $s_{0}=0$, then the relation (1-3) will be:\\\\
\begin{equation}
\sum_{n\leq x}a(n)+\frac{1}{2}a(x)=\frac{1}{2\pi
i}\int_{b-iT}^{b+iT}A(s)\frac{x^{s}}{s}ds+O(\frac{x^{b}B(b)}{T})+\\O(xH(2x)min(1,\frac{log(x)}{T}))\
\end{equation}
\textbf{Lemma 3.2}\\\\
\textit{Let $c_{1}$ be a positive constant, then $\zeta(s,\eta)$
has no zero point for $\sigma \geq
1-\frac{c_{1}}{log(|t|+2)}$ where $s=\sigma+it$.}\\\\
\textit{\textbf{Proof:}}\\\\ Let's define a new function as:
\begin{equation}
\zeta(s,\hat{\eta})=\sum_{n=1}^{\infty}\frac{(1-\eta(n))}{n^{s}}\
\end{equation}
If $\sigma \geq 1-\frac{c_{1}}{log(|t|+2)}$ then can be written:
\begin{equation}
|\zeta(s,\eta)|\leq |\zeta(s,\hat{\eta})|\
\end{equation}
Also,
\begin{equation}
|\zeta(s,\eta)+\zeta(s,\hat{\eta})|\leq |\zeta(s)|\
\end{equation}
Using Rouche's Theorem, the number of zeros of the function
$\zeta(s,\eta)$ is equal to the zeros of the function $\zeta(s)$.
Since $\zeta(s)$ has no zero in the set $\sigma \geq
1-\frac{c_{1}}{log(|t|+2)}$, therefore $\zeta(s,\eta)$ has also no
zero in the same set. This proves the Lemma.\\\\ \textbf{Lemma
3.3}\\\\ \textit{Let Dirichlet's series
$\Psi(s,\eta)=\sum_{n=1}^{\infty}\frac{\eta(n)\Lambda(n)}{n^{s}}$,
then, $\sigma_{c}= \sigma_{a}=1$ where $\sigma_{c}$ and
$\sigma_{a}$ denote the convergence and absolute convergence
abscissa series
respectively.}\\\\
\textit{\textbf{Proof:}}\\\\
Using the Prime Number Theorem for Arithmetic Progressions and
$gcd(k,N)=1$ one can write by referring to [11]:
\begin{equation}
\sum_{n\leq x,~n\equiv N
mod(k)}~\Lambda(n)=\frac{x}{\phi(k)}+O(\frac{x}{log^{H}x})\
\end{equation}
\begin{equation}
\lim_{x\rightarrow\infty}{\sum_{n\leq x,~n\equiv N
mod(k)}~\Lambda(n)}=\frac{x}{\phi(k)}\
\end{equation}
For sufficiently large number $N$:
\begin{equation}
\eta(n)=\prod_{p>\sqrt{N}, gcd(p,N)=1}(1-\frac{1}{\phi(p)})\
\end{equation}
Generally, one can write:
\begin{eqnarray}
\sum_{n\leq x}\eta(n)\Lambda(n)=\sum_{n\leq x}\prod_{p>\sqrt{N},
gcd(p,N)=1}(1-\frac{1}{\phi(p)})\Lambda(n)\
\end{eqnarray}\
and for sufficiently large number $N$ is;
\begin{eqnarray}
\sum_{n\leq x,~x\rightarrow\infty}\eta(n)\Lambda(n)=\sum_{n\leq
x,~x\rightarrow\infty}~~~\Lambda(n)\prod_{p>\sqrt{N},
gcd(p,N)=1}(1-\frac{1}{\phi(p)})=\nonumber\\ \lim_{x\rightarrow\infty}\psi(x)\left\{\prod_{p>\sqrt{N}~~,
gcd(p,N)=1}(1-\frac{1}{p-1})\right\}=~~~~~~~~~~~~~~~\nonumber\\
x\left\{\prod_{p>\sqrt{N}~~,
gcd(p,N)=1}(1-\frac{1}{p-1})\right\}~~~~~~~~~~~~~~~~~\
\end{eqnarray}\
Just,we give a conjecture namely the Sabihi's second
conjecture on the Goldbach's conjecture as below (to see our first
conjecture refer to [8]):
\\\\ \textbf{Sabihi's
Second Conjecture (SSC)}\\\\\textit{Let $\gamma$ be the Euler
constant,$p$ a prime number,$N$ a sufficiently large even
number,and Riemann Hypothesis holds,then:}
\begin{equation}
log(N)=\frac{4e^{-\gamma}\prod_{p>2}(1-\frac{1}{(p-1)^2})\prod_{p>2~,p\mid
N}\frac{p-1}{p-2}(1+O(\frac{1}{log(N)}))}{\prod_{p>\sqrt{N}, gcd(p,N)=1}(1-\frac{1}{p-1})}\
\end{equation}
Having the conjecture the relation (10-3) gives:
\begin{eqnarray}
\sum_{n\leq x}\eta(n)\Lambda(n)=4e^{-\gamma}\times~~~~~~~~~~~~~~~~~~~~~~\nonumber\\
\prod_{p>2}(1-\frac{1}{(p-1)^2})\prod_{p>2~,p\mid
N}\frac{p-1}{p-2}(1+O(\frac{1}{log(N)}))\frac{x}{log(N)}\
\end{eqnarray}
Assume we hold $N$ to a constant value and $x\rightarrow \infty$
then could be written:
\begin{equation}
\sigma_{c}=\lim_{x\rightarrow \infty}\frac{log|\sum_{n\leq
x}\eta(n)\Lambda(n)|}{log(x)}=1\
\end{equation}\
Therefore, $\sigma_{c}=\sigma_{a}=1$ and the Lemma is
proven.\\\\
\textbf{Lemma 3.4}\\\\ \textit{Let Dirichlet's series in
Lemma 3.3 have a pole at $s=1$,then it has a residue of the following form at same pole:}
\begin{equation}
4e^{-\gamma}\prod_{p>2}(1-\frac{1}{(p-1)^2})\prod_{p>2~,p\mid
N}\frac{p-1}{p-2}(1+O(\frac{1}{log(N)}))\frac{1}{log(N)}\
\end{equation}\
\textit{\textbf{Proof:}}\\\\
It is well-known
\begin{equation}
\lim_{s\rightarrow
1}(1-s)\sum_{n=1}^{\infty}\frac{\Lambda(n)}{n^{s}}=1\
\end{equation}\
On the other hand
\begin{equation}
\lim_{x\rightarrow\infty}\sum_{n\leq x,~n\equiv N
mod(k)}\frac{\Lambda(n)}{n^{s}}=\frac{1}{\phi(k)}\sum_{n=1}^{\infty}\frac{\Lambda(n)}{n^{s}}
\end{equation}\
Consequently
\begin{equation}
\lim_{x\rightarrow\infty}\sum_{n\leq x,~n\equiv N
mod(k)}\frac{\eta(n)\Lambda(n)}{n^{s}}=\sum_{n=1}^{\infty}~~\prod_{p>\sqrt{N},~gcd(p,N)=1}(1-\frac{1}{p-1})\sum_{n=1}^{\infty}\frac{\Lambda(n)}{n^{s}}
\end{equation}\
and this is equal to
\begin{equation}
4e^{-\gamma}\prod_{p>2}(1-\frac{1}{(p-1)^2})\prod_{p>2~,p\mid
N}\frac{p-1}{p-2}(1+O(\frac{1}{log(N)}))\frac{1}{log(N)}\sum_{n=1}^{\infty}\frac{\Lambda(n)}{n^{s}}\
\end{equation}\
Just, one could obtain the residue of the Dirichlet's series as
below:
\begin{eqnarray}
\lim_{s\rightarrow
1}(1-s)\sum_{n=1}^{\infty}\frac{\eta(n)\Lambda(n)}{n^{s}}=4e^{-\gamma}\prod_{p>2}(1-\frac{1}{(p-1)^2})\prod_{p>2~,p\mid
N}\frac{p-1}{p-2}(1+O(\frac{1}{log(N)}))\times \nonumber\\
\frac{1}{log(N)}(\lim_{s\rightarrow
1}(1-s))\sum_{n=1}^{\infty}\frac{\Lambda(n)}{n^{s}}=4e^{-\gamma}\prod_{p>2}(1-\frac{1}{(p-1)^2})\times ~~~~~~~~~~~~~~~~\nonumber\\
\prod_{p>2~,p\mid
N}\frac{p-1}{p-2}(1+O(\frac{1}{log(N)}))\frac{1}{log(N)}~~~~~~~~~~~~~~~~~~~~~~~~~~~~~
\end{eqnarray}
\section{PROOF OF INEQUALITY $\pi(N,\eta)\geq 1$ FOR SUFFICIENTLY LARGE EVEN INTEGER $N$}
In this section, we prove one main theorem. On
the basis of this theorem,the inequality $\pi(N,\eta)\geq 1$ under two
conditions will be proved.The first codition is to be assumed
trueness of Riemann Hypothesis (RH) and the second is to be
assumed to hold true Sabihi's second conjecture (SSC).\\\\
\textbf{Sabihi's Theorem}\\\\
\textit{Let $N$ be a sufficiently large even integer.Let 
both Sabihi's Second Conjecture (SSC) and Riemann Hypothesis (RH)
hold then:}
\begin{eqnarray}
\psi(N,\eta)=4e^{-\gamma}\prod_{p>2}(1-\frac{1}{(p-1)^2})\prod_{p>2~,p\mid
N}\frac{p-1}{p-2}(1+O(\frac{1}{log(N)}))\frac{N}{log(N)}+\nonumber\\O(N\{H(2N)min(1,\frac{log(N)}{T})-e^{-c\sqrt{log(N)}}\})~~~~~~~~~~~~~~~~~~~~~~~\
\end{eqnarray}\
\begin{eqnarray}
\pi(N,\eta)=4e^{-\gamma}\prod_{p>2}(1-\frac{1}{(p-1)^2})\prod_{p>2~,p\mid
N}\frac{p-1}{p-2}(1+O(\frac{1}{log(N)}))\frac{N}{log^{2}N}+~~~~~~~~~~~~~~~~~~~~~~~~~~~~~~~~~~~~~~~~~~~~~\nonumber\\\int_{2}^{N}\frac{du}{log^{2}u}+O(\frac{xH(2N)min(1,\frac{log(N)}{T})-\frac{Ne^{-c\sqrt{log(N)}}}{log(N)}-\frac{N^\frac{1}{2}}{log(N)}}{log(N)})+~~~~~~~~~~~~~~~~~~~~~~~~~~~~~~~~~~~~~~~~~~~~~~~~\nonumber\\O(\int_{2}^{N}\frac{(e^{-c\sqrt{log(u)}}-u^\frac{-1}{2})du}{log^{2}u})~~~~~~~~~~~~~~~~~~~~~~~~~~~~~~~~~~~~~~~~~~~~~~~~~~~~~~~~~~~~~~~~~~~~~~~~~~\nonumber\\
\end{eqnarray}\
\textit{\textbf{Proof:}}\\\\
Let $a=1-c_{1}\frac{1}{log(T+2)},
b=1+\frac{1}{log(x)},logT=(log(x)^{\frac{1}{\alpha+1}}$ for
$0<\alpha<1$,and $H(u)\leq logu$, $B(u)\leq c_{2} log(x)$ where
$c_{2}$ is a positive constant.By applying the Lemma 3.1,the
relation (2-3)and assuming $A(s)=\sum_{n=1}^{\infty}a(n)n^{-s}$
from Lemma 3.1:
\begin{equation}
\sum_{n\leq x}a(n)=\sum_{n\leq x}\eta(n)\Lambda(n)=\psi(x,\eta)\
\end{equation}\
Let $C_{n}$ be a rectangle contour with vertices $a\pm iT$ and
$b\pm iT$ then if $T$ tends to infinity
\begin{equation}
I=\frac{1}{2\pi i}\int_{C_{n}}\Psi(s,\eta)\frac{x^{s}}{s}ds=J+K\
\end{equation}\
where $J$ denotes integral along the line joining $b-iT$ to $b+iT$
and K denotes the integral along the other three sides of
rectangle.Applying the relation (2-3)to the relation (23-4) gives
\begin{equation}
I=\frac{1}{2\pi
i}\int_{b-iT}^{b+iT}\Psi(s,\eta)\frac{x^{s}}{s}ds+O(f(x))+O(g(x))=
\end{equation}\
\begin{equation}
J=\frac{1}{2\pi i}\int_{b-iT}^{b+iT}\Psi(s,\eta)\frac{x^{s}}{s}ds\
\end{equation}\
and
\begin{equation}
K=\frac{1}{2\pi
i}\{\int_{b-iT}^{a-iT}+\int_{a-iT}^{a+iT}+\int_{a+iT}^{b+iT}\}\Psi(s,\eta)\frac{x^{s}}{s}ds\
\end{equation}\
\begin{equation}
O(f(x))=\frac{xlog(x)}{T}_{T\rightarrow \infty}=0\
\end{equation}\
\begin{eqnarray}
O(g(x))=(\int_{b-i\infty}^{a-i\infty}+\int_{a-i\infty}^{a+i\infty}+\int_{a+i\infty}^{b+i\infty})\Psi(s,\eta)\frac{x^{s}}{s}ds\leq \nonumber\\
(\int_{b-i\infty}^{a-i\infty}+\int_{a-i\infty}^{a+i\infty}+\int_{a+i\infty}^{b+i\infty})\sum_{n=1}^{\infty}\Lambda(n)n^{-s}\frac{x^{s}}{s}ds=\nonumber\\
(\int_{b-i\infty}^{a-i\infty}+\int_{a-i\infty}^{a+i\infty}+\int_{a+i\infty}^{b+i\infty})(Z(s))
\frac{x^{s}}{s}ds\nonumber\\
\ll(\int_{b-i\infty}^{a-i\infty}+\int_{a-i\infty}^{a+i\infty}+\int_{a+i\infty}^{b+i\infty})log^{2}|t|\frac{x^{s}}{s}ds\
\end{eqnarray}\
Referring to [10],one can conclude that $z(s)\ll log^{2}|t|$ since
$\sigma>1-c_{1}log^{2}|t|$. By manipulating the last right hand
term of inequality(28-4) we obtain as below:
\begin{equation}
O(g(x))=xe^{-c\sqrt{log(x)}}\
\end{equation}\
Just, residue theorem expresses that
\begin{equation}
I=4e^{-\gamma}\prod_{p>2}(1-\frac{1}{(p-1)^2})\prod_{p>2~,p\mid
N}\frac{p-1}{p-2}(1+O(\frac{1}{log(N)}))\frac{x}{log(N)}\
\end{equation}\
The relation (2-3) gives:
\begin{equation}
\sum_{n\leq
x}a(n)=\psi(x,\eta)=I-O(f(x))-O(g(x))+O(\frac{x^{b}B(b)}{T})+O(xH(2x)min(1,\frac{log(x)}{T}))\
\end{equation}\
Consequently
\begin{eqnarray}
\psi(x,\eta)=4e^{-\gamma}\prod_{p>2}(1-\frac{1}{(p-1)^2})\prod_{p>2~,p\mid
N}\frac{p-1}{p-2}(1+O(\frac{1}{log(N)}))\frac{x}{log(N)}+\nonumber\\O(x\{H(2x)min(1,\frac{log(x)}{T})-e^{-c\sqrt{log(x)}}\})~~~~~~~~~~~~~~~~~~~~~~~\
\end{eqnarray}\
If we apply $x=N$ for a sufficiently large even integer, then the
first formula of the theorem is proven. By applying and combining
the relations(12-2)to(14-2)and(32-4)could be written:
\begin{eqnarray}
\pi(x,\eta)=\frac{\psi(x,\eta)-O(x^\frac{1}{2})}{log(x)}+\int_{2}^{x}\frac{\psi(u,\eta)-O(u^\frac{1}{2})}{ulog^{2}(u)}du=4e^{-\gamma}\prod_{p>2}(1-\frac{1}{(p-1)^2})\times~~~~~~~~~~~~~~~~~~~~~~~~~~~~\nonumber\\
\prod_{p>2~,p\mid
N}\frac{p-1}{p-2}(1+O(\frac{1}{log(N)}))\frac{x}{log(N)log(x)}+~~~~~~~~~~~~~~~~~~~~~~~~~~~~~~~~~~~~~~~~~~~~~~~~~\nonumber\\O(\frac{xH(2x)min(1,\frac{log(x)}{T})}{log(x)})-O(\frac{xe^{-c\sqrt{log(x)}}}{log(x)})+\int_{2}^{x}\frac{du}{log^{2}u}+O(\int_{2}^{x}\frac{e^{-c\sqrt{log(u)}}du}{log^{2}u})-~~~~~~~~~~~~~~~~~~~~~~~~~~~~~\nonumber\\O(\frac{x^\frac{1}{2}}{log(x)})-O(\int_{2}^{x}\frac{du}{u^\frac{1}{2}log^{2}u})~~~~~~~~~~~~~~~~~~~~~~~~~~~~~~~~~~~~~~~~~~~~~~~~~~~~~~~~~~~\
\end{eqnarray}\
Again, applying $x=N$ to the relation (33-4) for a sufficiently
large even integer, the second formula of the theorem is also
proven. In the above relation, the following relation and
inequality have been applied.Applying the Prime Number Theorem:
\begin{equation}
\psi(x)=x+O(xe^{-c\sqrt{log(x)}})\
\end{equation}
and
\begin{equation}
\int_{2}^{x}\frac{\psi(u,\eta)du}{ulog^{2}u}\leq
\int_{2}^{x}\frac{\psi(u)du}{ulog^{2}u}\
\end{equation}
From the relation (33-4) and replacing $x$ by $N$ and tending $N$ to infinity, we easily see
that $\pi(N,\eta)>1$ or $\pi(N,\eta)\neq 0$ and the theorem is
proven.This proves the strong Goldbach's conjecture for
sufficiently large even integers.

\end{document}